\documentclass[preprint,12pt,times]{article}
\usepackage{graphicx}
\usepackage{amsthm}
\usepackage{amssymb}
\usepackage{amsmath}
\usepackage{amscd}
\usepackage{bbm}
\usepackage{float}
\usepackage{booktabs}
\numberwithin{equation}{section}

\newtheorem{theorem}{Theorem}[section]
\newtheorem{lemma}{Lemma}[section]

\usepackage[numbers,sort&compress]{natbib}

\numberwithin{figure}{section}
\numberwithin{table}{section}

\allowdisplaybreaks[2]
\newcommand{\ep}{\mathrm{e}}

\newcommand\btd{\raise 2pt \hbox{$\hat\bigtriangledown$}\hskip 1.5pt}
\newcommand\bt{\raise 2pt \hbox{$\bigtriangledown$}\hskip 1.5pt}

\begin{document}
\date{}
\title{ Blow-up analysis for the generalized hyperelastic rod equation }

\author{Shaojie Yang\thanks{Corresponding author: shaojieyang@kust.edu.cn (Shaojie Yang)} \\~\\
\small Department of  Mathematics,~~Kunming University of Science and Technology,  \\
\small Kunming, Yunnan 650500, China\\~\\
\small Research Center for Mathematics and Interdisciplinary Sciences,\\
 \small Kunming University of Science and Technology,\\
\small Kunming, Yunnan 650500, China}

\date{}
\maketitle
\begin{abstract}
~\\We investigate the blow-up  for the generalized hyperelastic rod equation, which includes the Camassa-Holm equation, the the hyperelastic-rod wave equation,  and the rotation-Camassa-Holm equation as special cases.
We establish  a new local-in-space blow-up criterion (i.e., a criterion involving only the properties of the initial data $u_0$ in a small neighborhood of a single point $x_0\in\mathbb{R}$) under certain initial data. Our local-in-space blow-up condition on the initial data is purely local in the space variable.\\


\end{abstract}
\noindent\rule{15.5cm}{0.5pt}

\section{Introduction}\label{sec1}
In this paper, we consider the following generalized hyperelastic rod equation
\begin{align}\label{a1}
u_t-u_{txx}+\left[f(u)\right]_x-\left[f(u)\right]_{xxx}+\left[g(u)+\frac{f''(u)}{2}u_x^2\right]_x=0,
\end{align}
where $f,g:\mathbb{R}\rightarrow\mathbb{R}$ are two given $C^\infty$-functions.  Eq.\eqref{a1} was proposed by Holden and Raynaud in \cite{r4}. Tian, Yan and Zhang \cite{r3} studied the local well-posedness and  blow-up  for \eqref{a1}.  Brandolese and Cortez \cite{r2} considered the  local-in-space blow-up issues for \eqref{a1}.  The global conservative weak solutions for \eqref{a1} after wave breaking was shown in \cite{r5}. The uniqueness of conservative   solution was proved in \cite{r33}.

For $f(u)=\frac{u^2}{2}$ and $g(u)=u^2$, Eq.\eqref{a1} becomes the  celebrated  Camassa-Holm (CH) equation 
\begin{align}\label{CH}
u_t-u_{txx}+3uu_x=2u_xu_{xx}+uu_{xxx},
\end{align}
which is a  shallow water waves model describing the unidirectional propagation. The CH equation was theoretically derived formally by Fokas and Fuchssteiner in \cite{r8}, but later was physically derived from the Euler’s  equations  by Camassa and Holm \cite{r1} in 1993.  It admits a famous feature of peaked solitons (peakons). The CH equation has been widely investigated  extensively in the last three decades from qualitative standpoint, such as well-posedness \cite{r9,r10,r11}, persistence properties \cite{r18,r19},  blow-up of solutions \cite{r14,r15,r16,r17,Rr1}, orbital stability \cite{r12,r13,R30}, global conservative and dissipative solutions \cite{r20,r21}.

For $f(u)=\frac{\gamma u^2}{2}$ and $g(u)=\frac{3-\gamma}{2}u^2$, Eq.\eqref{a1} becomes the hyper-elastic rod wave equation 
\begin{align}\label{hrq}
u_t-u_{txx}+3uu_x=\gamma(2u_xu_{xx}+uu_{xxx}),
\end{align}
which was introduced by Dai \cite{r6,r7}, and  describes the propagation of nonlinear waves inside cylindrical hyperelastic rods and  $u(t,x)$ represents the radial stretch relative to a pre-stressed state.
The well-posedness and  blow-up criteria for \eqref{hrq} have been investigated in \cite{r22, R22}. 
Local-in-space  blow-up criteria for \eqref{hrq} has been shown in \cite{r23}. The  orbital stability of solitary waves for Eq.\eqref{hrq} has been proved in \cite{r24}.

For $f(u)=\frac{u^2}{2}+\frac{\beta_0}{\beta}u$ and $g(u)=(c-\frac{\beta_0}{\beta})u+u^2+\frac{\omega_1}{3\alpha^2}u^3+\frac{\omega_2}{4\alpha^3}u^4$,  Eq.\eqref{a1} becomes the rotation-Camassa-Holm (R-CH) equation
\begin{align}\label{rch}
u_t-u_{txx}+cu_x+3uu_x+\frac{\omega_1}{\alpha^2}u^2u_x+\frac{\omega_2}{\alpha^3}u^3u_x-\frac{\beta_0}{\beta}u_{xxx}=2uu_{xxx}+u_xu_{xx},
\end{align}
which was derived from incompressible and irrotational two-dimensional equatorial shallow water with the Coriolis effect under the CH scaling \cite{r25}.
Global existence and uniqueness of the energy conservative weak solutions  to the R-CH equation have been derived in \cite{r26}.  Well-posedness, travelling waves and geometrical aspects  have been studied in \cite{r27}.
Non-uniform dependence and well-posedness in the sense of Hadamard have been proved in \cite{r28}.
Generic regularity of conservative solutions has been studied in \cite{r29}.  Wang, Yang and Han \cite{r30} proved that symmetric waves to the R-CH equation must be traveling waves.  Existence and stability of solitary waves to the R-CH equation
have been shown in \cite{r31}.

In Ref.\cite{r23,R1}, Brandolese   introduced a local-in-space criteria for blow-up in the study of the hyperelastic rod wave equation \eqref{hrq},  this local-in-space blow-up criteria  only involves the values of $u'_0 (x_0)$ and $u_0 (x_0) $ in a single point $x_0$. 
Later, Brandolese and Cortez \cite{r2} establish a local-in-space blowup criterion for the generalized hyper-elastic rod equation \eqref{a1}, the main result  is as follows.
\begin{theorem}[\cite{r2}]\label{T1}
Let $u_0\in H^s(\mathbb{R})$ with $s>3/2$. Let $f,g\in C^\infty(\mathbb{R})$ with $f''\geq \gamma>0.$ The maximal time
$T^*$ of the solution $u$ to  Eq.\eqref{a1}  in $C([0,T^*);H^s(\mathbb{R}))\cap C^1([0,T^*);H^{s-1}(\mathbb{R}))$ must be finite, if at least one of the two following conditions {\rm (1)} or {\rm (2)} is fulfilled:\\~\\
 \noindent{\rm(1)}
\begin{minipage}[t]{0.9\linewidth}
$\exists $ $c\in\mathbb{R}$ such that $m=g(c)=\min\limits_{\mathbb{R}} g$.\\
 The map $\phi: \mathbb{R} \to \mathbb{R}$ given by $\phi=\sqrt{\frac{g-m}{\gamma}}$ is $K$-Lipschitz with $0\leq K\leq 1$.\\
 $\exists~x_0\in \mathbb{R}$ such that
 \begin{align}
  u_0'(x_0)<-\frac{1}{2}\left(\sqrt{1+8K^2}-1 \right)|u_0(x_0)-c|.
  \end{align}
  \end{minipage}\\
   \noindent{\rm(2)} 
\begin{minipage}[t]{0.9\linewidth}
Or, otherwise,\\
$\exists $ $c\in\mathbb{R}$ such that $M=g(c)=\max\limits_{\mathbb{R}} g$.\\
 The map $ \psi: \mathbb{R} \to \mathbb{R}$ given by $\psi=\sqrt{\frac{M-g}{\gamma}}$ is $K$-Lipschitz with $0\leq K\leq \frac{1}{\sqrt{8}}$.\\
 $\exists~x_0\in \mathbb{R}$ such that
 \begin{align}
  u_0'(x_0)<-\frac{1}{2}\left(1-\sqrt{1-8K^2} \right)|u_0(x_0)-c|.
  \end{align}
  \end{minipage} \\~\\
 More precisely, the following upper bound estimate for  $T^*$ holds:
  \begin{align*}
  T^*\leq \frac{4}{\gamma\sqrt{4u_0'(x_0)^2-(\sqrt{1\pm8K^2}-1)^2(u_0(x_0)-c)^2           }},
  \end{align*}
 where in the term $\pm 8K^2$ one has to take the positive sign under the conditions of Part {\rm (1)} and the
negative sign under the conditions of Part {\rm (2)}.
\end{theorem}

In this paper, we  establish a new local-in-space blow-up for \eqref{a1},  which is an extension and improvement of the previous works by Brandolese and Cortez \cite{r2}. The main result of this paper is stated as follows.
\begin{theorem}\label{T}
Let $u_0\in H^s(\mathbb{R})$ with $s>3/2$. Let $f,g\in C^\infty(\mathbb{R})$ with $f''\geq \gamma>0.$ The maximal time
$T^*$ of the solution $u$ to the Cauchy problem \eqref{a2} in  $C([0,T^*);H^s(\mathbb{R}))\cap C^1([0,T^*);H^{s-1}(\mathbb{R}))$ must be finite, if at least one of the two following conditions {\rm (i)} or {\rm (ii)} is fulfilled:\\~\\
 \noindent{\rm(i)}
\begin{minipage}[t]{0.9\linewidth}
$\exists $ $c\in\mathbb{R}$ such that $m=g(c)=\min\limits_{\mathbb{R}} g$.\\
 The map $\phi: \mathbb{R} \to \mathbb{R}$ given by $\phi=\sqrt{\frac{g-m}{\gamma}}$ is $K$-Lipschitz with $0\leq K\leq 1$.\\
 $\exists~x_0\in \mathbb{R}$ such that
 \begin{align}
  u_0'(x_0)<-\frac{1}{2K}\left(\sqrt{1+8K^2}-1 \right)\phi(u_0(x_0)).
  \end{align}
    More precisely, the following upper bound estimate for  $T^*$ holds:
  \begin{align*}
  T^*\leq \frac{4K}{\gamma\sqrt{4K^2u_0'(x_0)^2-(\sqrt{1+8K^2}-1)^2\phi(u_0(x_0))^2           }}.
  \end{align*}    
   \end{minipage}\\
   \noindent{\rm(ii)} 
\begin{minipage}[t]{0.9\linewidth}
Or, otherwise,\\
$\exists $ $c\in\mathbb{R}$ such that $M=g(c)=\max\limits_{\mathbb{R}} g$.\\
 The map $ \psi: \mathbb{R} \to \mathbb{R}$ given by $\psi=\sqrt{\frac{M-g}{\gamma}}$ is $K$-Lipschitz with $0\leq K\leq \frac{1}{\sqrt{8}}$.\\
 $\exists~x_0\in \mathbb{R}$ such that
 \begin{align}
  u_0'(x_0)<-\frac{1}{2K}\left(1-\sqrt{1-8K^2} \right)\psi(u_0(x_0)).
  \end{align}
  More precisely, the following upper bound estimate for  $T^*$ holds:
  \begin{align*}
  T^*\leq \frac{4K}{\gamma\sqrt{4K^2u_0'(x_0)^2-(\sqrt{1- 8K^2}-1)^2\psi(u_0(x_0))^2 .       }}.
  \end{align*}     \end{minipage} \\~\\
\end{theorem}

\noindent {\bf Remark.~}
Note that the $K$-Lipschitz Lipschitz condition on $\phi$ and $\psi$, that is 
\begin{align*}
|\phi(u_0(x_0))|=|\phi(u_0(x_0))-\phi(c)|\leq K|u_0(x_0)-c|
\end{align*}
and
\begin{align*}
|\psi(u_0(x_0))|=|\psi(u_0(x_0))-\psi(c)|\leq K|u_0(x_0)-c|,
\end{align*}
which implies that
\begin{align*}
-\frac{1}{2}\left(\sqrt{1+8K^2}-1 \right)|u_0(x_0)-c|\leq -\frac{1}{2K}\left(\sqrt{1+8K^2}-1 \right)\phi(u_0(x_0))\end{align*}
and 
\begin{align*}
  -\frac{1}{2}\left(1-\sqrt{1-8K^2} \right)|u_0(x_0)-c|\leq -\frac{1}{2K}\left(1-\sqrt{1-8K^2} \right)\psi(u_0(x_0)).
  \end{align*}
  Consequently, we enlarge blow-up upper bound of  $ u_0'(x_0)$, then  our  blow-up condition covers blow-up condition by  Brandolese and Cortez in  \cite{r2}.

The remainder of this paper is organized as follows. In Section \ref{sec2}, we recall several useful results which are crucial in  in the proof of  Theorem \ref{T}. In Section 3, we prove Theorem \ref{T}.

\section{Preliminaries}\label{sec2}

In this section, we recall several useful results which are crucial in  in the proof of  Theorem \ref{T}.

First of all, using $p:=\frac{1}{2}\ep^{-|x|}$, it is convenient to write the Cauchy problem associated with
the Eq.\eqref{a1} in the following nonlocal form
\begin{equation}\label{a2}
\begin{cases}
u_t+f'(u)u_x+\partial_xp*\left[g(u)+\frac{f''(u)}{2}u_x^2\right]=0,\\
u(0,x)=u_0(x).
\end{cases}
\end{equation}

The local well-posedness for the Cauchy problem \eqref{a2} is as follows.
\begin{lemma}[\cite{r3}]\label{le1}
Assume that $f,g\in C^{\infty}(\mathbb{R})$ and  $u_0(x)\in H^s(\mathbb{R})(s>3/2)$, then there exist $T=T(u_0,f,g)>0$ and a unique solution $u\in C([0,T);H^s(\mathbb{R}))\cap C^1([0,T);H^{s-1}(\mathbb{R}))$ of the Cauchy problem \eqref{a2}. The solution has constant energy integral
\begin{align*}
\int_\mathbb{R}(u^2+u_x^2)dx=\int_\mathbb{R}(u_0^2+u_{0x}^2)dx=\|u_0\|_{H^1}^2.
\end{align*}
Moreover, the solution depends continuously on the initial data, i.e., the mapping $u_0\rightarrow u(\cdot,u_0):H^s\rightarrow C([0,T);H^s(\mathbb{R}))\cap C^1([0,T);H^{s-1}(\mathbb{R}))$ is continuous.
\end{lemma}

Next, the necessary and sufficient condition for the blowup  to the Cauchy problem \eqref{a2} can be formulated as follows.
\begin{theorem}[\cite{r3}]\label{le2}
Assume that $f,g\in C^{\infty}(\mathbb{R})$ and  $f''\geq\gamma>0$. Given $u_0(x)\in H^s(\mathbb{R})(s>3/2)$, then the solution $u$ of the Cauchy problem \eqref{a2} blows up in finite time $T^*<+\infty$ if and only if
\begin{align*}
\liminf_{t\to T^*} \inf_{x\in \mathbb{R}} \{ u_x(t,x) \}=-\infty.
\end{align*}
\end{theorem}

Finally, we recall the following convolution estimates  which will be used in the sequel.
\begin{lemma}[\cite{r2}]\label{L3}
Let $1_{\mathbb{R}^\pm}$ denote one of the two indicator functions $1_{\mathbb{R}^+}$ or $1_{\mathbb{R}-}$ .\\~\\
   \noindent{\rm(1)} 
\begin{minipage}[t]{0.9\linewidth}
If $f$ and $g$ satisfies the condition as in Theorem \ref{T1}{\rm (i)} then the following estimate holds:
\begin{align}
\left(p1_{\mathbb{R}^\pm} \right)*\left( g(u)+\frac{f''(u)}{2}u_x^2\right)\geq \frac{\alpha}{2}(g(u)-m)+\frac{m}{2}
\end{align}
with
\begin{align}\label{b1}
\alpha=\frac{1}{4K^2}(\sqrt{1+8K^2}-1).
\end{align}
\end{minipage}\\~\\
 \noindent{\rm(2)} 
\begin{minipage}[t]{0.9\linewidth}
If $f$ and $g$ satisfies the condition as in Theorem \ref{T1}{\rm (ii)} then the following estimate holds:
\begin{align}
\left(p1_{\mathbb{R}^\pm} \right)*\left( g(u)+\frac{f''(u)}{2}u_x^2\right)\geq \frac{\alpha}{2}(g(u)-M)+\frac{M}{2}
\end{align}
with
\begin{align*}
\alpha=\frac{1}{4K^2}(1-\sqrt{1-8K^2}).
\end{align*}
\end{minipage}\\~\\
In the case $g = m = M$ be a constant function (this corresponds to $K = 0$), the right-hand side
of the above convolution estimates reads 
\begin{align*}
\left(p1_{\mathbb{R}^\pm} \right)*\left( g(u)+\frac{f''(u)}{2}u_x^2\right)\geq \frac{g(u)}{2}.
\end{align*}
\end{lemma}

\section{Blow-up}\label{sec3}

In this section, we prove our local-in-space blow-up Theorem \ref{T}.\\

\noindent\textbf{Proof of Theorem \ref{T}}.  Taking the space derivative $x$ in \eqref{a2}, we obtain
\begin{align}
u_{tx}+f'(u)u_{xx}=-\frac{f''(u)}{2}u_x^2+g(u)-p*\left[g(u)+\frac{f''(u)}{2}u_x^2 \right].
\end{align}
We consider the following flow map
\begin{align}
\begin{cases}
\frac{q(t,x)}{dt}=f'(u(t,q(t,x))),~~&x\in\mathbb{R}, t>0,\\
q(0,x)=x,~~&x\in\mathbb{R}
\end{cases}
\end{align}
Notice that the assumptions made on $f$ and $u$ imply that $q\in C^1([0,T^*)\times\mathbb{R}, \mathbb{R})$  is well defined on the whole time interval $[0, T^*).$

We first proceed putting the conditions of Theorem \ref{T}(i). Note that the uniform convexity condition $f''\geq\gamma>0$. By the definition of  $\alpha$ in \eqref{b1} yields $0<\alpha\leq 1$, then we have
\begin{align}\label{b2}
\notag\frac{d}{dt}[u_x(t,q(t,x))]=&[u_{tx}+f'(u)u_{xx}](t,q(t,x))\\
\notag=&-\frac{f''(u)}{2}u_x^2+g(u)-p*\left( g(u)+\frac{f''(u)}{2}u_x^2\right)\\
\notag\leq&\left[(1-\alpha)(g(u)-m)-\frac{\gamma}{2}u_x^2 \right](t,q(t,x))\\
=&\frac{\gamma}{2}\left( 4K^2\alpha^2\phi^2(u)-u_x^2\right)(t,q(t,x)),
\end{align}
where  $\phi=\sqrt{\frac{g(u)-m}{\gamma}}$.

Let 
\begin{align*}
A(t,x)=\left(2K\alpha \phi(u)-u_x\right)(t,q(t,x))
\end{align*}
and
\begin{align*}
B(t,x)=\left(2K\alpha \phi(u)+u_x\right)(t,q(t,x)),
\end{align*}
then \eqref{b2} implies that
\begin{align*}
\notag\frac{d}{dt}[u_x(t,q(t,x))]\leq \frac{\gamma}{2}(AB)(t,x).
\end{align*}
Realling that
\begin{align*}
p=p1_{\mathbb{R}^+}+p1_{\mathbb{R}^-},~~~~p_x=p1_{\mathbb{R}^-}-p1_{\mathbb{R}^+}
\end{align*}
and the inequality
\begin{align*}
f''\geq \gamma,
\end{align*}
then
\begin{align*}
A_t(t,x)=&2K\alpha \phi'(u)(u_t+f'(u)u_x)-(u_{tx}+f'(u)u_{xx})\\
=&\frac{f''(u)}{2} u_x^2-g(u)+\left(p-2K\alpha\phi'(u)p_x \right)*\left( g(u)+\frac{f''(u)}{2}u_x^2\right)\\
\geq &\frac{\gamma}{2}u_x^2-g(u)+\left(1+2K\alpha\phi'(u)\right)p1_{\mathbb{R}^+}*\left( g(u)+\frac{f''(u)}{2}u_x^2\right)\\
&+\left(1-2K\alpha\phi'(u)\right)p1_{\mathbb{R}^-}*\left( g(u)+\frac{f''(u)}{2}u_x^2\right)\\
\geq& \frac{\gamma}{2}u_x^2+(\alpha-1)(g(u)-m)\\
=&-\frac{\gamma}{2}\left( 4K^2\alpha^2\phi^2(u)-u_x^2\right)\\
=&-\frac{\gamma}{2}AB,
\end{align*}
and 
\begin{align*}
B_t(t,x)=&2K\alpha \phi'(u)(u_t+f'(u)u_x)+(u_{tx}+f'(u)u_{xx})\\
=&-\frac{f''(u)}{2} u_x^2+g(u)+\left(p+2K\alpha\phi'(u)p_x \right)*\left( g(u)+\frac{f''(u)}{2}u_x^2\right)\\
\geq &-\frac{\gamma}{2}u_x^2+g(u)-\left(1+2K\alpha\phi'(u)\right)p1_{\mathbb{R}^-}*\left( g(u)+\frac{f''(u)}{2}u_x^2\right)\\
&-\left(1-2K\alpha\phi'(u)\right)p1_{\mathbb{R}^+}*\left( g(u)+\frac{f''(u)}{2}u_x^2\right)\\
\geq& -\frac{\gamma}{2}u_x^2-(\alpha-1)(g(u)-m)\\
=&\frac{\gamma}{2}\left( 4K^2\alpha^2\phi^2(u)-u_x^2\right)\\
=&\frac{\gamma}{2}AB,
\end{align*}
where we used the convolution estimates in Lemma \ref{L3} and $\phi$-is $K$-Lipschitz with $0\leq K\leq1$.

According to the assumption on initial data, that is 
\begin{align*}
u_0'(x_0)<-\frac{1}{2K}(\sqrt{1+8K^2}-1)\phi(u_0(x_0)),
\end{align*}
it's easy to find that
\begin{align*}
A(0, x_0)>0~~{\rm and}~~B(0,x_0)<0.
\end{align*}
By the continuity of $A(t,x_0)$ and $B(t,x_0)$, then it ensures that
\begin{align*}
A_t(t, x_0)>0~~{\rm and}~~B_t(t,x_0)<0,~~\forall~~t\in [0,T),
\end{align*}
which implies that
\begin{align*}
A(t,x_0)>A(0,x_0)>0~~{\rm and}~~B(t,x_0)<B(0,x_0)<0, ~~\forall~~t\in [0,T).
\end{align*}

Let
\begin{align*}
h(t)=\sqrt{-(AB)(t,x_0)},
\end{align*}
a simple calculation yields
\begin{align*}
\frac{dh(t)}{dt}=&\frac{-A_tB-AB_t}{2h}(t,x_0)\\
\geq& \frac{\gamma(-AB)(A-B)}{4h}(t,x_0)\\
\geq&\frac{\gamma}{2}h^2(t), 
\end{align*}
where we used the inequality $\frac{A-B}{2}\geq h$. Hence, the solution blows up infinite time $T^*$ with 
\begin{align*}
T^*\leq \frac{2}{\gamma\sqrt{u_0'^2(x_0)-4K^2\alpha^2\phi^2(u_0(x_0)) }}.
\end{align*}

The condition for Theorem \ref{T} (ii) is obtained in a very similar way to  above, we omit it here. Hence,  the proof
of Theorem \ref{T} is completed.

\section*{Acknowledgment}
This work is supported by the National Natural Science Foundation of China (Grant No. 12561101).

\section*{Conflict of interest }
The author do not have any other competing interests to declare.

\section*{Data availability }
No data was used for the research described in the article.

\end{document}